\documentclass[11pt,draft]{amsart}

\usepackage{amssymb}
\usepackage{url}
\usepackage{amscd}
\usepackage{amsthm}
\usepackage{amsmath}

\theoremstyle{plain}

\theoremstyle{remark}

\numberwithin{equation}{section}
\begin{document}

\date{}

\title[Collinear triples in permutations]
{Collinear triples in permutations   }
\author{Liangpan Li}
\date{}
\address{Department of Mathematics,
Shanghai Jiaotong University, Shanghai 200240, People's Republic of
China} \email{liliangpan@yahoo.com.cn}

\begin{abstract}
Let $\alpha:\mathbb{F}_q\rightarrow\mathbb{F}_q$ be a permutation
and $\Psi(\alpha)$ be the number of collinear triples in the graph
of $\alpha$, where $\mathbb{F}_q$ denotes a finite field of $q$
elements. When $q$ is odd Cooper and Solymosi once proved
$\Psi(\alpha)\geq(q-1)/4$ and conjectured the sharp bound should be
$\Psi(\alpha)\geq(q-1)/2$. In this note we indicate that the
Cooper-Solymosi conjecture is true.
\end{abstract}

\maketitle

\section{The main result and its proof}
Let $\mathbb{F}_q$ be the finite field of $q$ elements with $q$ odd.
Let $\alpha:\mathbb{F}_q\rightarrow\mathbb{F}_q$ be a permutation
and $\Psi(\alpha)$ be the number of collinear triples in
\[G_{\alpha}=\{(i,\alpha(i)):i\in\mathbb{F}_q\},\] the
graph of $\alpha$. Cooper and Solymosi \cite{CooperSolymosi} first
(see also \cite{li}) obtained the lower bound
\begin{equation}\label{0.25}
\Psi(\alpha)\geq\frac{q-1}{4},
\end{equation}
and conjectured the best one should be
\begin{equation}\label{0.5}
\Psi(\alpha)\geq\frac{q-1}{2}.
\end{equation}
Later Cooper \cite{Hypergraphs} showed that the problem of counting
collinear triples in a permutation and the finite plane Kakeya
problem are intimately connected, and improved (\ref{0.25}) into
\begin{equation}
\Psi(\alpha)\geq\frac{5q-1}{14}.
\end{equation}
This is the right way to solve the Cooper-Solymosi conjecture and we
shall explain Cooper's idea in more detail in Section 3.

A subset in $\mathbb{F}_q^2$ containing a line in each direction is
called a Kakeya set. According to Cooper, given the permutation
$\alpha$, one can construct a corresponding Kakeya set
\[K_{\alpha}=\bigcup_{i\in\mathbb{F}_{q}}L(i,(0,\alpha(i))),\] where
$L(s,x)$ denotes the line in $\mathbb{F}^2_q$ through $x$ with slope
$s$. Let $\Gamma_{\alpha}$ be the hypergraph on the vertex set
$G_{\alpha}$ whose edges are the maximal collinear subsets of
$G_{\alpha}$, and write $\|\Gamma_{\alpha}\|$ for the quantity
\[\sum_{e\in E(\Gamma_{\alpha})}{|e|-1\choose 2}.\]
Then by applying an incidence formula of Faber \cite{faber}, Cooper
successfully showed that
\begin{equation}\label{formula}
\sharp K_{\alpha}=\frac{q(q+1)}{2}+\|\Gamma_{\alpha}\|.
\end{equation}
To confirm (\ref{0.5}), by considering
\[\Psi(\alpha)=\sum_{e\in E(\Gamma_{\alpha})}{|e|\choose 3}\geq\sum_{e\in E(\Gamma_{\alpha})}{|e|-1\choose 2}=\|\Gamma_{\alpha}\|,\]
it suffices to prove \[\sharp
K_{\alpha}\geq\frac{q(q+1)}{2}+\frac{q-1}{2}.\] So we are led to
caring about the size of the finite plane Kakeya sets.

The finite field Kakeya problem posed by Wolff in his influential
survey \cite{wolff} asks for the smallest size a Kakeya set can
have. Coincidentally, for any Kakeya set $K\subset\mathbb{F}_q^2$
Faber \cite{faber} once proved the bound
\[\sharp K\geq\frac{q(q+1)}{2}+\frac{q}{3},\]
and conjectured the sharp one should be
\begin{equation}\label{kakeya}\sharp
K\geq\frac{q(q+1)}{2}+\frac{q-1}{2}.\end{equation}  Recently, this
problem was solved by Blokhuis and Mazzocca \cite{blokhuis}, see
also \cite{Hypergraphs} for other medium bound. As an immediate
corollary, the Cooper-Solymosi conjecture (\ref{0.5}) turns out to
be true.

To make this note  more self-contained, the incidence formula of
Faber (a new proof) and the connection between the permutations and
the finite plane Kakeya sets discovered by Cooper (with more
details) will be discussed in the next two sections. We refer the
reader \cite{blokhuis} for the original proof of the Faber
conjecture (\ref{kakeya}).

\section{The Incidence formula of Faber}

Suppose $K$ is a minimal Kakeya set in $\mathbb{F}^2_q$. Clearly, we
may assume $K$ is of the form
\[K=\bigcup_{s\in PG(1,q)}L(s,x^{(s)}).\]
 Let $\mu_{x}$ be the number of these lines passing
through $x\in\mathbb{F}^2_q$. Obviously one has
\[
\mu_{x}=\sum_{s\in PG(1,q)}\chi_{L_{s}}(x),\] where  we let
$\chi_{A}$ denote the characteristic function of
$A\subset\mathbb{F}^2_q$ and write $L_{s} =L(s,x^{(s)})$ for
simplicity. The incidence formula of Faber \cite{faber} says that
\begin{equation}\label{faber formula}\sharp K=\frac{q(q+1)}{2}+\sum_{x\in K}{\mu_{x}-1\choose
2}.\end{equation} In the following we will give a succinct proof of
the Faber formula. As we know, two lines with different slopes
intersect at one point. Hence considering
\begin{align}
\sum_{x\in K}\mu_{x}&=\sum_{x\in K}\sum_{s\in
PG(1,q)}\chi_{L_{s}}(x)\label{step 1}\\ &=\sum_{s\in
PG(1,q)}\sum_{x\in
K}\chi_{L_{s}}(x)\nonumber\\&=q(q+1),\nonumber\\
\sum_{x\in K}\mu_{x}^2&=\sum_{x\in K}\sum_{i\in PG(1,q)}\sum_{j\in
PG(1,q)}\chi_{L_{i}}(x)\chi_{L_{j}}(x)\label{step 2}\\&=\sum_{i\in
PG(1,q)}\sum_{j\in
PG(1,q)}\sum_{x\in K}\chi_{L_{i}}(x)\chi_{L_{j}}(x)\nonumber\\
&=\sum_{i\in PG(1,q)}\sum_{j\in PG(1,q)}\sharp(L_i\cap
L_j)\nonumber\\&=2q(q+1),\nonumber\end{align} it follows that
\begin{align*}\sum_{x\in
K}\frac{(\mu_x-1)(\mu_x-2)}{2}&=\frac{\sum_{x\in
K}\mu_{x}^2}{2}-\frac{3\sum_{x\in K}\mu_{x}}{2}+\sharp K\\
&=\frac{2q(q+1)}{2}-\frac{3q(q+1)}{2}+\sharp K\\
&=\sharp K-\frac{q(q+1)}{2}.\end{align*} Finally we indicate that
(\ref{step 1}) and (\ref{step 2}) already appeared in
\cite{iosevich}.

\section{The collinear tuple hypergraphs and  the finite plane Kakeya sets }

As before, given the permutation $\alpha$, one can construct a
corresponding Kakeya set
\[K_{\alpha}=\bigcup_{i\in\mathbb{F}_{q}}L(i,(0,\alpha(i))).\]
Let $x=(x_1,x_2)\in\mathbb{F}_{q}^2$ be any point satisfying
$\mu_x\geq3.$ Since $\alpha$ is a permutation, $x_1>0$. We assume
that $x$ lies in the lines
\[L(i_k,(0,\alpha(i_k)))\ \ (k=1,2,\ldots,\mu_x),\]
which means
\[\frac{x_2-\alpha(i_k)}{x_1-0}=i_k.\]
Hence for all $1\leq j<k\leq \mu_{x}$,
\[\frac{\alpha(i_k)-\alpha(i_j)}{i_k-i_j}=-x_1.\]
In another words, the set
\[E_x\doteq\big\{(i_k,\alpha(i_k))\big\}_{k=1}^{\mu_x}\]
is collinear. In fact, \[E_x\subset L(-x_1,(i_1,\alpha(i_1)))\] is a
maximal collinear subset of $G_{\alpha}$. For if
\[\frac{\alpha(k)-\alpha(i_1)}{k-i_1}=-x_1\]
holds for some $k\neq i_1$, then
\[\alpha(k)=\alpha(i_1)-(k-i_1)x_1=x_2-kx_1,\]
which means $x$ lies in the line $L(k,(0,\alpha(k)))$. Thus $k=i_j$
holds for some $2\leq j\leq \mu_x$ and $E_x$ is a maximal collinear
subset of $G_{\alpha}$.

On the other hand, suppose $\{(i_{t},\alpha(i_{t})\}_{t=1}^{\gamma}$
is a maximal collinear subset of $G_{\alpha}$, where $\gamma\geq3$.
Define
\[z_1=\frac{\alpha(i_2)-\alpha(i_1)}{i_1-i_2}\]
and
\[z_2=\alpha(i_1)+z_1i_1,\]
then it is easy to verify that
\[(z_1,z_2)\in L(i_t,\alpha(i_t))\]
holds for $t=1,2,\ldots,\gamma$.

In summary, the point $x$ in $K_{\alpha}$ with $\mu_x\geq3$
corresponds to a collinear $\mu_{x}$-tuple in $G_{\alpha}$, and vice
visa. Consequently, \begin{equation}\label{kakeyahyper}\sum_{x\in
K_{\alpha}}{\mu_{x}-1\choose 2}=\sum_{e\in
E(\Gamma_{\alpha})}{|e|-1\choose 2}.\end{equation} Combining
(\ref{kakeyahyper}) with the Faber formula (\ref{faber formula})
yields (\ref{formula}).

\section{Acknowledgements}
The author thanks Aart Blokhuis, Qing Xiang and Yaokun Wu for kindly
pointing out the recent progresses on the finite field Kakeya
problem to him.

\end{document}